\documentclass[10pt,letterpaper,oneside,journal]{IEEEtran}

\usepackage{amsmath}
\usepackage{amsthm}
\usepackage{amssymb}
\usepackage{graphicx}
\usepackage{cases}
\usepackage{balance}
\usepackage{empheq}  

\usepackage{here}

\usepackage{enumitem}
\newcommand{\er}{\varepsilon}

\usepackage{xcolor}
\definecolor{magenca}{rgb}{0,8,0.0,8} 

\theoremstyle{plain}
\newtheorem{theorem}{Theorem}

\newtheorem{proposition}{Proposition}
\newtheorem{assumption}{Assumption}
\newtheorem{lemma}{Lemma}
\newtheorem{remark}{Remark}

\newtheorem{assumption*}{Assumption}

\usepackage[english]{babel}

\DeclareMathOperator\sign{sgn}
\newcommand{\scal}[2]{\left\langle #1,\, #2\right\rangle  }
\newcommand{\norm}[1]{\left\Vert#1 \right\Vert}

\begin{document}

\title{Boundary Feedback and Observer Synthesis for a Class of Nonlinear Parabolic--Elliptic PDE Systems}
\author{ Kamal Fenza, Moussa Labbadi and Mohamed Ouzahra
\thanks{K. Fenza and M. Ouzahra are with the  Laboratory of  Mathematics and Applications to Engineering Sciences,  Sidi Mohamed Ben Abdellah University, Fes, Morocco}
\thanks{M. Labbadi is with the Aix-Marseille University, LIS UMR CNRS 7020, 13013 Marseille, France.}
}

\maketitle
\IEEEoverridecommandlockouts

\maketitle \thispagestyle{empty} 
 
\maketitle

\begin{abstract}
This paper investigates the stabilization of a coupled system comprising a parabolic PDE and an elliptic PDE with  nonlinear terms. A rigorous backstepping design provides an explicit boundary control law and exponentially convergent observers from partial boundary measurements. Several theorems ensure exponential stability and well-posedness of the nonlinear closed-loop system.
 \end{abstract}
\vspace{-1cm}
\section{Introduction}
\vspace{-0.5cm} Nonlinear parabolic–elliptic systems appear frequently in a wide variety of physical and biological contexts \cite{Ahn2019, Wu2006}. The nonlinear couplings inherent in these systems often generate complex dynamics that pose significant challenges to both theoretical analysis and control design \cite{TaoWinkler2015, Parada2020}. Unlike linear systems, nonlinearities can give rise to instabilities and intricate behaviors, demanding sophisticated mathematical tools for their study. In this note, we extend the results from \cite{morris2025boundary} to the nonlinear case and tackle the problems of both control and observer design within this more general framework. Related contributions in the literature include \cite{KrsticSmyshlyaev2008, Parada2020, VazquezKrstic2009, VazquezTrelatCoron2008}.

The first contribution of this work is the design of a single boundary feedback control law that ensures exponential stabilization of the coupled nonlinear system. This contrasts with previous approaches that often rely on multiple control inputs. By employing a backstepping transformation originally developed for parabolic PDEs (e.g., \cite{KrsticSmyshlyaev2008}), we derive an explicit expression for the stabilizing controller. The kernel properties of the transformation facilitate the construction of a target system in parabolic–elliptic form, bypassing the need for eigenfunction expansions. Under specific inequalities involving the system parameters, exponential stability of the closed-loop system is rigorously established. The second contribution is the design of an observer that estimates the full state of the system using a single boundary measurement. The observer is constructed via a backstepping approach tailored to the coupled nonlinear setting. The exponential convergence of the observation error is proven under suitable conditions, again leveraging a transformation adapted from the parabolic PDE context. Finally, the closed-loop stability of the nonlinear system under output feedback is thoroughly analyzed.

Building on the framework introduced in \cite{morris2025boundary}, we consider the following class of nonlinear parabolic--elliptic systems:
\begin{gather}
u_t(x, t) = u_{xx}(x, t) - \rho u(x, t)+ f_1( u(x, t))\nonumber\\
\quad  +\alpha v(x, t)+f_2( v(x, t)),\label{eq:eq1}\\
0  = v_{xx}(x, t) - \gamma v(x, t) + \beta u(x, t)+f_3( u(x, t)), \label{eq:eq2} \\
    u_x(0, t) = 0,  u_x(1, t) = \omega(t), 
    v_x(0, t) = 0,  v_x(1, t) = 0, 
\label{eq:eq3} \\
    u(x, 0) = u_0(x), \quad v(x, 0) = v_0(x),
\label{eq:eq4} 
\end{gather}
for \( x \in [0,1] \) and \( t \geq 0 \). Where, \( \rho, \alpha, \beta, \gamma \in \mathbb{R} \) are system parameters, while the nonlinear mappings \( f_1, f_2, f_3 : L^2(0,1) \to L^2(0,1) \) represent distributed nonlinear functions, which can also be interpreted as perturbations in the coupled dynamics. Equation~\eqref{eq:eq1} defines a parabolic PDE, whereas equation~\eqref{eq:eq2} is elliptic, characterized by the absence of time dynamics but the presence of second-order spatial derivatives. The system coupling arises whenever \( \alpha \) or \( \beta \) are nonzero or the nonlinearities \( f_2, f_3 \) are active. Control is exerted via a Neumann boundary input \( \omega(t) \) at the boundary \( x=1 \) of the parabolic subsystem. In our analysis, we consider consistent initial conditions, meaning that \( u_0 \), \( v_0 \)  and $f_3(u_0)$ must satisfy the elliptic equation  
\begin{equation*}
    0 = (v_0)_{xx} - \gamma v_0 + \beta u_0 + f_3(u_0).
\end{equation*}
Failure to satisfy this condition results in impulsive behavior at \( t = 0 \), which is not physically meaningful.\\
Our objectives are twofold: to design a boundary feedback control law that guarantees exponential stabilization of the coupled nonlinear system, and to develop an observer relying solely on the boundary measurement \( u(1,t) \). This observer facilitates accurate state estimation despite limited measurement access, which is critical for output feedback stabilization.

In the absence of nonlinearities (i.e., \( f_1 = f_2 = f_3 = 0 \)), the system reduces to the linear parabolic--elliptic model studied in \cite{morris2025boundary}, whereboundary control and observer design were thoroughly developed. The present work extends these results to incorporate nonlinearities, which introduce significant challenges in control synthesis and estimation.


\section{Main results}
\subsection{Boundary control  for (\ref{eq:eq1})--(\ref{eq:eq4})}
In this section, we present the main results of the paper, including the control and observer design, along with the analysis of their stability and well-posedness. To mitigate this, we design a single stabilizing boundary control input $w(t)$. Unlike the approaches used in the linear case including~\cite{KrsticSmyshlyaev2008} (see Chapter~10) and~\cite{VazquezKrstic2006, VazquezKrstic2009}, which employed two  boundary actuators to stabilize coupled parabolic--elliptic PDEs without nonlinear disturbances, we achieve stabilization of the coupled system~(\ref{eq:eq1})--(\ref{eq:eq4}).
We  begin with well-posedness and stability of  the open-loop system (\ref{eq:eq1})--(\ref{eq:eq4}). Let  $\Delta z(x) = \frac{d^2 z}{dx^2}$.	We define the operator $A_{\gamma} : D(A_{\gamma} ) \rightarrow L^2(0, 1)$ with domain $D(A_{\gamma} ) = \left\{ z \in H^2(0, 1) \ \bigg| \ z'(0) = z'(1) = 0 \right\}$, by $A_{\gamma}  z = (\gamma I - \Delta)z,$ where$\gamma \neq -(n\pi)^2,\quad n = 0, 1, 2, \ldots$, for all $z \in D(A_{\gamma} ) $. Thus the linear  operator $(A_{\gamma})^{-1} : L^2(0, 1) \rightarrow D(A_{\gamma} )$, is bounded. 	Alternatively,  we define  the operator $A$ by 
	\begin{equation}\label{11X}
		A = \Delta - \rho I + \alpha \beta(A_{\gamma})^{-1},
	\end{equation}
	with domain
	$$
	D(A) = \left\{ w \in H^2(0, 1) \ \bigg| \ w'(0) = w'(1) = 0 \right\}.
	$$      
  The operator $A$  generates the $C_0$-semigroup $T(t)$ on $L^2(0,1)$.  Furthermore, according to (\cite{morris2025boundary}, p. 3) , $T(t)$ is exponentially stable if and only if $\rho > \frac{\alpha \beta}{\gamma}$.  
    Let us consider the following assumption: 
\begin{assumption}\label{assuption:Lischitz_F}
     We assume that $f_1,f_2,f_2$ are globally Lipschitz functions such that $f_1(0)=f_2(0)=f_3(0)=0$, with $M_1$, $M_2$, and $M_3$ be the Lipschitz constants of $f_1$, $f_2$ and $f_3$, respectively.
\end{assumption}
We present the first result as follows:
\begin{theorem}\label{loop-sys}
    	Under Assumption~\ref{assuption:Lischitz_F},  $\gamma \neq -(n\pi)^2 $ \text{with}  $n\in \mathbb{N}$ and  If
        \begin{enumerate}
            \item  The semigroup $T(t)$ is exponentially stable.
            \item $M=\rho-\frac{\alpha \beta}{\gamma }-M_1 +M_3 | \alpha|+  M_2 |\beta|+ M_2M_3 >0 $.
        \end{enumerate}
		Then, the system \eqref{eq:eq1}--\eqref{eq:eq4} with \( w(t) = 0 \) is exponentially stable, admitting the following upper bound:
			\[
			\| u(\cdot, t) \| + \| v(\cdot, t) \| \leq ce^{-M t}\left( \| u_0 \| + \| v_0 \| \right) ,
			\]
			where$c$ is a positive constant independent of $M$, $\|u_0\|$ and $\|v_0\|$.
\end{theorem}
When the semigroup $T(t)$ is not exponentially stable, we build a feedback control of system \eqref{eq:eq1}--\eqref{eq:eq4}. Our approach relies on an invertible state transformation, previously employed in backstepping designs, to convert the original system into a target system with desirable stability properties.

\begin{equation}\label{KAMAL}
				\tilde{u}(x, t) = u(x, t) - \int_0^x k(x, y)\, u(y, t)\, dy
		\end{equation}
		on the parabolic state $u(x, t)$, while the elliptic state $v(x, t)$ is unchanged and  the kernel $k(x, y)$ is given  in Lemma \ref{lem:SmyshlyaevKrstic2004}.\\
Following Lemma \ref{Lemma3}, we define the inverse transformation of (\ref{KAMAL}) as:
\begin{equation}\label{Ivertible}
u(x, t) = \tilde{u}(x, t) + \int_0^x \ell(x, y) \tilde{u}(y, t) \, dy,
\end{equation}
\begin{theorem}\label{thm:target_system}
Let $k(x, y)$   be given by the kernel  in Lemma~\ref{lem:SmyshlyaevKrstic2004}, and consider the following  control signal:
\begin{equation}\label{control}
\omega(t) = \int_0^1 k_x(1, y) u(y, t) \, dy + k(1, 1) u(1, t). 
\end{equation}
Then transformation (\ref{KAMAL}) converts the system (\ref{eq:eq1})--(\ref{eq:eq4}) into the following target system:
\begin{gather}  
	 \tilde{u}_t(x,t)=\tilde{u}_{xx}(x, t) - (c_1 + \rho)\, \tilde{u}(x, t) + \alpha v(x, t)\notag\\
			 - \alpha \int_0^x k(x, y)\, v(y, t)\, dy+	 f_1\left(\tilde{u}(x, t)+\int_0^x  \ell(x, y)\, \tilde{u}(y, t)\, dy\right)\notag\\+f_2(v(x, t))- \int_0^x k(x, y)f_2(v(y, t))dy  \notag\\  -\int_0^x k(x, y)  f_1\left(\tilde{u}(y, t)+\int_0^y  \ell(y, z)\, \tilde{u}(y, t)\, dz\right)dy,   \notag\\ 	
0 = v_{xx}(x, t) - \gamma v(x, t) + \beta\, \tilde{u}(x, t) + \beta \int_0^x  \ell(x, y)\, \tilde{u}(y, t)\, dy+\notag\\f_3\left( \tilde{u}(x, t) + \int_0^x  \ell(x, y)\, \tilde{u}(y, t)\, dy\right),\notag\\
			\tilde{u}_x(0, t) = 0, \ \tilde{u}_x(1, t) = 0,   \tilde{u}(x,0)=\tilde{u}_0 \
v(x,0)=v_0. \label{eq:eq755}
	\end{gather}
\end{theorem}
The purpose of the following theorem is to study the well-posedness of the coupled system~(\ref{eq:eq1})--(\ref{eq:eq4})  withe the input (\ref{control}). 
\begin{theorem}\label{thm:well-posdnes}
Under Assumption~\ref{assuption:Lischitz_F},   $\gamma \neq -(n\pi)^2 $ \text{with}  $n\in \mathbb{N}$. Then, the target system (\ref{eq:eq755}) admits a unique mild solution $(\tilde{u}, v)\in C([0,+\infty), L^2(0, 1)) \times C([0,+\infty), L^2(0, 1))$.
    \end{theorem}
Next, we provide conditions that ensure the exponential stability of the target system~\eqref{eq:eq755}.
\begin{theorem}\label{thm:stability}
    Let $ u_0, v_0 \in L^2(0, 1)$ be consistent initial conditions. Assume that Assumption~\ref{assuption:Lischitz_F} and   $\gamma \neq -(n\pi)^2 $ \text{with}  $n\in \mathbb{N}$. If  $K_1=c_1 + \rho- \big(M_1+(M_2+|\alpha |)(|\beta| +M_3 )  \big) (1 +N_{c_1} )^2  >0 $, then system (\ref{eq:eq755}) is exponentially stable.
		Furthermore,   the exponential decay is specified by
			\(
			\| \tilde{u}(\cdot, t) \| + \| v(\cdot, t) \| \leq c \left( \| \tilde{u}_0 \| + \| v_0 \| \right) e^{-K_1 t},
			\)
			where $c$ is a positive constant independent of $K_1$, $\|\tilde{u}_0\|$ and $\|v_0\|$.
\end{theorem}
The main result of this section is now immediate.
	\begin{theorem}\label{Theo of sys1}
Let  \( u_0, v_0 \in L^2(0, 1) \) and let \( k(x, y) \) be given by Lemma \ref{lem:SmyshlyaevKrstic2004}. If the parameter  $K_1 >0$  and  Assumption~\ref{assuption:Lischitz_F} hold, then the control signal~ (\ref{control}) exponentially stabilizes the closed-loop system (\ref{eq:eq1})--(\ref{eq:eq4}) in the space \( L^2(0,1) \times L^2(0,1) \), with the convergence rate \(
\|u(\cdot, t)\| + \|v(\cdot, t)\| \leq c e^{-K_1 t} \left( \|u_0\| + \|v_0\| \right),  \forall t\geq 0.\)
and \( c \) is a positive constant independent of \( K_1 \), \( \|u_0\| \),  \( \|v_0\| \).
\end{theorem}
\vspace{-0.5cm}
	\subsection{State observer for (\ref{eq:eq1})--(\ref{eq:eq4})}
    Since the control law \eqref{control} requires full-state information, a state observer is necessary. In distributed parameter systems, domain-wide measurements are rarely available—only boundary measurements are typically accessible. We design a backstepping-based observer for \eqref{eq:eq1}–\eqref{eq:eq4}, using only the boundary measurement \( u(1,t) \) to reconstruct the state.
We propose the following observer for system (\ref{eq:eq1})--(\ref{eq:eq4}):
    \begin{gather*}
		\hat{u}_t(x, t) = \hat{u}_{xx}(x, t) - \rho \hat{u}(x, t) + f_1(\hat{u}(x, t)) \\
	 + \alpha \hat{v}(x, t) + f_2(\hat{v}(x, t)) 
		+ \sigma_1(x)[u(1, t) - \hat{u}(1, t)],  \\
		0 = \hat{v}_{xx}(x, t) - \gamma \hat{v}(x, t) + \beta \hat{u}(x, t)+ f_3(\hat{u}(x, t)),  \\
		\hat{u}_x(0, t) = 0, \quad \hat{u}_x(1, t) = \omega(t) + \sigma_2[u(1, t) - \hat{u}(1, t)], \\
		\hat{u}(x, 0) = \hat{u}_0(x), \quad \hat{v}(x, 0) = \hat{v}_0(x).
	\end{gather*}
	The in-domain output injection function $\sigma_1(\cdot)$  and  boundary injection value $\sigma_2$  are to be designed. 
    Define the error states
	\begin{gather*}
		\er^u(x, t) = u(x, t) - \hat{u}(x, t), \ \er^v(x, t) = v(x, t) - \hat{v}(x, t).
	\end{gather*}
 Then, the observation error satisfies   
\begin{gather}
		\er^u_t(x, t) = \er^u_{xx}(x, t) - \rho \er^u(x, t) + \alpha \er^v(x, t) - \sigma_1(x) \er^u(1, t) \nonumber\\
		\quad f_1(u(x, t))  -f_1\big(u(x, t)-\er^u(x, t)\big)+ f_2(v(x, t)) \nonumber\\ -f_2\big(v(x, t)-\er^v(x, t)\big), \label{eq:obs_err1} \\
		0 = \er^v_{xx}(x, t) - \gamma \er^v(x, t) + \beta \er^u(x, t) \nonumber \\
		\quad +f_3(u(x, t))  -f_3\big(u(x, t)-\er^u(x, t)\big),\\
\er^u_x(0, t) = 0, \quad \er^u_x(1, t) = -\sigma_2 \er^u(1, t) \\
\er^v_x(0, t) = 0, \quad \er^v_x(1, t) = 0\label{eq:obs_errn}, 
\end{gather}
    	with consistent initial conditions. A backstepping approach is used to select $\sigma_1(\cdot)$ and $\sigma_2$ so that the error system~\eqref{eq:obs_err1}--\eqref{eq:obs_errn} is exponentially stable.
We define the backstepping state transformation as follows:
	\begin{equation}\label{tran11}
	    \er^{\tilde{u}}(x, t) = \er^{u}(x, t) -\int_0^x k(x, y) \er^{u}(y, t) dy,
	\end{equation}
    where$k( . , .)$ is given by Lemma \ref{lem:SmyshlyaevKrstic2004}. By Lemma \ref{lem:morris2025boundary}, the
inverse transformation is
\begin{equation}\label{tran12}
	    \er^{u}(x, t) = \er^{\tilde{u}}(x, t) +\int_0^x \ell(x, y) \er^{\tilde{u}}(y, t) dy,
	\end{equation}

	\begin{proposition}\label{obsever_pro}
		If the output injections are
		\begin{equation}
			\sigma_1(x) = 0, \quad \sigma_2 = -k(1, 1), 
		\end{equation}
		 then transformation (\ref{tran11}) converts the error system~(\ref{eq:obs_err1})--(\ref{eq:obs_errn}) into the  following target system:
		\begin{gather}
		 \er^{\tilde{u}}_t(x, t) = \er^{\tilde{u}}_{xx}(x, t) - (p_1+ \rho) \er^{\tilde{u}}(x, t) + \alpha \er^{\tilde{v}}(x, t)\nonumber\\
		\quad -  \int_0^x k(x, y)  \left( f_1(u(y, t))  -f_1(u(y, t)-\er^{u}(y, t) ) \right) dy \nonumber\\
        -\int_0^x k(x, y)  \left( f_2(v(y, t))  -f_2(v(y, t)-\er^{v}(y, t) ) \right)dy
        \nonumber \\
   -f_2(v(x, t)-\er^{v}(x, t) )- \alpha \int_0^x k(x, y)\, \er^v(y, t)\, dy \nonumber \\
          +  f_2(v(x, t))+ f_1(u(x, t))  -f_1(u(x, t)-\er^{u}(x, t) ), \label{eq:o11}\\  
		0 = \er^{v}_{xx}(x, t) - \gamma \er^{v}(x, t) + \beta \er^{u}(x, t)+f_3(u(x, t)) \nonumber \\
		\quad -f_3\left(u(x, t)-\er^{u}(x, t) \right), \nonumber\\
\er^{\tilde{u}}_x(1, t) = -\int_0^1 k_x(1, y) \left( \er^{\tilde{u}}(y, t) - \int_0^y \ell(y, z) \er^{\tilde{u}}(z, t) dz \right) dy, \label{eq:o2} \\
\er^{\tilde{u}}_x(0, t) = 0,  \er^{\tilde{v}}_x(0, t) = 0, \quad \er^{\tilde{v}}_x(1, t) = 0\label{eq:o3}.
\end{gather}
	\end{proposition}
The following theorem establishes the exponential stability of the observation error system (12)-(15) :
\begin{theorem}\label{stability_Observer}  Let  the assumption (\ref{assuption:Lischitz_F}) holds and   $\gamma \neq -(n\pi)^2 $ \text{with}  $n\in \mathbb{N}$. \\
		If 
\begin{equation}\label{condition}
    K_3=c_1 + \rho- K_4 >0, 
\end{equation}        
with $K_4=\left(( M_2 +|\alpha |  ) (|\beta| +M_3 ) +\frac{\eta^2+1}{2}+ M_1\right) (1+N_{c_1})^2$. Then system (12)-(15)  is exponentially stable.\\
		Furthermore,   the exponential decay is specified by
			\[
			\| \er^{u}(\cdot, t) \| + \| \er^{v}(\cdot, t) \| \leq c e^{-K_3 t} \left( \| \er^{u}(\cdot, 0) \| + \| \er^{v}(\cdot, 0)  \| \right),
			\]
			where $c$ is a positive constant independent of $K_3$.
\end{theorem}
\subsection{Output feedback control laws}
In practice, the full state of a system is typically not accessible for feedback control. Output feedback strategies aim to stabilize the system using only the available measurements. A widely used method consists of combining a state feedback controller—designed under the assumption of full state availability—with an observer that reconstructs the unmeasured state variables. 
Using only the boundary measurement \( u(1,t) \), the resulting output feedback controller couples the observer of system \eqref{eq:eq1}--\eqref{eq:eq4} with the following boundary feedback law:
\begin{equation}\label{contW}
    w(t) = \int_0^1 k_x(1, y)\, \hat{u}(y, t)\, dy + k(1,1)\, \hat{u}(1, t),
\end{equation}
where\( k(x, y) \) is the kernel in Lemma \ref{lem:SmyshlyaevKrstic2004}. Given that the original system~\eqref{eq:eq1}--\eqref{eq:eq4} is well-posed.
\begin{theorem}
  For any consistent initial conditions $u_0, v_0 \in L^2(0,1)$, the closed-loop system consisting of \eqref{eq:eq1}--\eqref{eq:eq4}, together with the observer dynamics and control signal  (\ref{contW}) admits a unique solution    in $C((0, \infty),  L^2(0,1))$. Moreover, the origin is exponentially stable in $L^2(0,1)$.
\end{theorem}

\section{Proof of the main results}
\subsection{Proof of Theorem~\ref{loop-sys}}
\begin{proof}
\textbf{Step 1: Well-posedness of the  open-loop system:}
 For $\gamma \neq -(n\pi)^2 $,  $n\in \mathbb{N}$, we have that the solution to the elliptic equation (\ref{eq:eq2}) is:  
\begin{equation}\label{solutionv0}
    v(\cdot,t)= \beta (\gamma I - \Delta)^{-1}(u(\cdot, t))+ (\gamma I - \Delta)^{-1} f_3(u(\cdot, t)). 
\end{equation}
  When  $w(t)=0$, the system (\ref{eq:eq1})--(\ref{eq:eq4}) can be reformulated as 
\begin{gather}
\frac{d}{dt} u(x,t) = A u(x,t) + \Psi(u(x,t)), \quad (x, t) \in (0,1)\times\mathbb{R}^{+} \nonumber \\
u_x(0, t) = 0, \quad u_x(1, t) = 0, \quad t \in \mathbb{R}^{+} \label{SystLoop} \\
u(x,0) = u_0, \quad x \in (0,1) \nonumber \\
v(x,t) = \beta (\gamma I - \Delta)^{-1}(u(x, t)) + (\gamma I - \Delta)^{-1} f_3(u(x, t)), \nonumber
\end{gather}
where$\Psi(u(x,t))=  f_1(u(x, t))+\alpha (\gamma I - \Delta)^{-1} f_3(u(x, t)) +    f_2\bigg(\beta (\gamma I - \Delta)^{-1}(u(x, t)) + (\gamma I - \Delta)^{-1} f_3(u(x, t))\bigg)$.\\
To show that (\ref{SystLoop}) is well-posed, we'll show that the operator $\Psi$ is Lipschitz globally on $L^{2}(0,1)$. Let \( u_1, u_2 \in L^2(0,1) \). Then, we use the fact that $f_2$ is globally Lipschitz. This allows us to deduce that 
$$
\begin{aligned}
&| \Psi(u_1(x)) - \Psi(u_2(x))|\leq | f_1(u_1(x)) - f_1(u_2(x))|\\&+ | \alpha | |(\gamma I - \Delta)^{-1} \bigg(f_3(u_1(x)) -f_3(u_2(x))\bigg)|\\
&+M_2 |\beta (\gamma I - \Delta)^{-1}\bigg(u_1(x)-u_2(x)\bigg)|\\
&+M_2 |(\gamma I - \Delta)^{-1}\bigg(  f_3(u_1(x)) -f_3(u_2(x))\bigg)|. 
\end{aligned}
$$
Using the fact that $f_1$ and $f_3$ are global functions and that $(\gamma I - \Delta)^{-1}$ is a bounded operator, we can conclude that
$$
\begin{aligned}
&| \Psi(u_1(x)) - \Psi(u_2(x))|\leq M| u_1(x) - u_2(x)|, 
\end{aligned}
$$
where$M=M_1 +M_3 | \alpha|+  M_2 |\beta|+ M_2M_3 $. Then 
\begin{equation}\label{esti-psi}
    \| \Psi(u_1)-\Psi(u_2) \|\leq C_2\| u_1-u_2 \|,
\end{equation}
Therefore, the mapping \( u \mapsto \Psi \) is globally Lipschitz on \( L^2(0,1) \).\\
Hence, according to (\cite{pazy}, p. 184), we conclude that for all $\tilde{u}_0\in L^2(0,1)$, the system (\ref{parabolic}) has a unique mild solution $u\in C([0,+\infty), L^2(0, 1))$. \\
\textbf{Step 2: Exponential stability result:}
	From (\cite{morris2025boundary}, p. 2), the spectrum of the  operator $A$ (given   in \ref{11X} ) is given
by the simple eigenvalues 
$$
\lambda_n = -\rho + \frac{\alpha \beta}{\gamma + (n\pi)^2} - (n\pi)^2, \quad n = 0, 1, 2, \ldots
$$
and the corresponding eigenfunctions  $\phi_n(x)=cos(n\pi x)$ for $n\geq 0$. Moreover,   the operator $A$  generates a $C_0-$ semigroup which is defined as follows:
	$$
	T(t)=\sum_{n=0}^{\infty} e^{\lambda_n t} \langle \cdot, \phi_n\rangle \phi_n, 
	$$
	Let $t\geq 0$,   the mild solution  is given by 
	$$
	u(\cdot, t)=T(t)u_0+\int_0^t T(t-s)\Psi  (u(\cdot, t))ds, 
	$$
	We have that
		$$\| u(\cdot, t)\| \leq e^{-t(\rho-\frac{\alpha \beta}{\gamma })  } \| u_0\|+ \int_0^t e^{-(t-s)(\rho-\frac{\alpha \beta}{\gamma })  }\| \Psi   (u((\cdot, s))\|ds.$$
		Then 
		$$
		e^{t(\rho-\frac{\alpha \beta}{\gamma })}\| u(\cdot, t)\| \leq \| u_0\| + \int_0^t e^{s(\rho-\frac{\alpha \beta}{\gamma })}\| \Psi(u(\cdot, s))\|ds.
		$$
		Under uniquality (\ref{esti-psi}) and  using the fact  that $(\gamma I - \Delta)^{-1}$ is a linear bounded operator, 
		$$
		e^{t(\rho-\frac{\alpha \beta}{\gamma })}\| u(\cdot, t)\| \leq \| u_0\| +M \int_0^t e^{s(\rho-\frac{\alpha \beta}{\gamma })}\| u(\cdot, s)\|ds,
		$$
	From Gronwall’s inequality, we obtain
		$$
		e^{t(\rho-\frac{\alpha \beta}{\gamma })}\| u(\cdot, t)\| \leq e^{Mt} \| u_0\|.
		$$
		Then 
		$$
		\| u(\cdot, t)\| \leq  e^{t(M-(\rho-\frac{\alpha \beta}{\gamma }))} \| u_0\|.
		$$
		Hence 
		\begin{equation}\label{kam}
			\| u(\cdot, t)\| \leq  e^{-t(\rho-\frac{\alpha \beta}{\gamma }-M)} \| u_0\|.
		\end{equation}
Moreover, the space \( L^2(0,1) \times L^2(0,1) \) is equipped with the norm defined by:
        \[
	\left\|
	\begin{bmatrix}
		\tilde{u}(\cdot, t) \\
		v(\cdot, t)
	\end{bmatrix}
	\right\|_{L^2 \times L^2}
	=
	\sqrt{ \| \tilde{u}(\cdot, t) \|^2 + \| v(\cdot, t) \|^2 }.
	\]
On the other hand, we have that
\begin{equation}\label{v}
    \| v(\cdot, t) \|\leq \bigg( |\beta| + M_3\bigg) \| u(\cdot, t) \|.
\end{equation}
Therefore, from equations (\ref{v}) and  (\ref{kam}), we have
		\[
		\|u(\cdot, t)\|+ \|v(\cdot, t)\|
		\leq \sqrt{2}\bigg(1+|\beta| + M_3\bigg)e^{-t(\rho-\frac{\alpha \beta}{\gamma }-M)}  \bigg(\| u_0\|+ \| v_0\|  \bigg).
		\]
	This completes the proof.
    \end{proof}
\subsection{Proof of Theorem~\ref{thm:target_system}}
\begin{proof}
We differentiate (\ref{Ivertible}) with respect to $x$ twice and with respect to $t$:
    \begin{gather}
u_{xx}(x, t) = \tilde{u}_{xx}(x, t) + \int_0^x k_{xx}(x, y)u(y, t)\,dy  \nonumber\\
 + k_x(x, x)u(x, t) + \frac{d}{dx}k(x, x)u(x, t) + k(x, x)u_x(x, t)\label{vvv}
\end{gather}
and
\begin{align}\label{t}
&u_t(x, t) = \tilde{u}_t(x, t) + k(x, x)u_x(x, t) - k_y(x, x)u(x, t)\notag \\ &+ k_y(x, 0)u(0, t) 
+ \int_0^x \left( k_{yy}(x, y) - \rho\,k(x, y)\right)u(y, t)\,dy \notag \\&+ \alpha\int_0^x k(x, y)v(y, t)\,dy+\int_0^x k(x, y)f_2(v(y, t))dy \notag \\&+\int_0^x k(x, y)  f_1\left(\tilde{u}(y, t)+\int_0^y  \ell(y, z)\, \tilde{u}(y, t)\, dz\right)dy.
\end{align}
Here, 
$$
k_x(x, x) = \left.\frac{\partial}{\partial x}k(x, y)\right|_{x = y}, \quad
k_y(x, x) = \left.\frac{\partial}{\partial y}k(x, y)\right|_{x = y}, 
$$
and $\frac{d}{dx}k(x, x) = k_x(x, x) + k_y(x, x).$\\   
Substituting equations \eqref{vvv} and (\ref{t}) into (\ref{eq:eq1}), and using the fact that $k_y(x, 0) = 0$, we obtain:
\[
\begin{aligned} 
&\tilde{u}_t(x, t) = \tilde{u}_{xx}(x, t) - (c_1 + \rho)\,\tilde{u}(x, t) + \alpha v(x, t) \\
&\quad - \alpha \int_0^x k(x, y)v(y, t)\,dy 
+ \left(2\,\frac{d}{dx}k(x, x) + c_1\right)u(x, t) \\
&\quad + \int_0^x \left[k_{xx}(x, y) - k_{yy}(x, y) - c_1 k(x, y)\right]u(y, t)\,dy\\
&+	 f_1\left(\tilde{u}(x, t))+\int_0^x  \ell(x, y)\, \tilde{u}(y, t)\, dy\right)+f_2(v(x, t))\\
 &-\int_0^x k(x, y)  f_1\left(\tilde{u}(y, t)+\int_0^y  \ell(y, z)\, \tilde{u}(y, t)\, dz\right)dy\\
 &-\int_0^x k(x, y)f_2(v(y, t))dy.
\end{aligned}
\]
Furthermore, $k(x, y)$ is given by Lemma \ref{lem:SmyshlyaevKrstic2004}, then
\begin{gather*}
\tilde{u}_t(x, t) = \tilde{u}_{xx}(x, t) - (c_1 + \rho)\,\tilde{u}(x, t) + \alpha v(x, t) \\
- \alpha \int_0^x k(x, y)v(y, t)\,dy +	f_2(v(x, t)) \\
+ f_1\left(\tilde{u}(x, t) + \int_0^x \ell(x, y)\, \tilde{u}(y, t)\,dy\right)\\
 -\int_0^x k(x, y)  f_1\left(\tilde{u}(y, t)+\int_0^y  \ell(y, z)\, \tilde{u}(y, t)\, dz\right)dy\\
 -\int_0^x k(x, y)f_2(v(y, t))dy.
\end{gather*}
Moreover, we have that
\[
\tilde{u}_x(0,t)= u_x(0, t) - k(0, 0)u(0, t) = 0,
\]
and
\begin{gather*}
u_x(1, t) = u_x(1, t) - \int_0^1 k_x(1, y)u(y, t)\,dy - k(1, 1)u(1, t) \\
= w(t) - \int_0^1 k_x(1, y)u(y, t)\,dy - k(1, 1)u(1, t) 
= 0.
\end{gather*}
Finally, owing to the consistency of the initial conditions $u_0$ and $v_0$, transformation (\ref{Ivertible}) ensures that the corresponding transformed initial state $\tilde{u}_0$ and $v_0$ remain consistent as well.
\end{proof}
\subsection{Proof of Theorem~\ref{thm:well-posdnes}}  
   \begin{proof}   
    \textbf{Step 1: Explicit form of the solution to the elliptic system:}
With the notation $\Delta z(x) = \frac{d^2 z}{dx^2}$,	we define the operator $A^\gamma : D(A^\gamma) \rightarrow L^2(0, 1)$ with domain $D(A^\gamma) = \left\{ z \in H^2(0, 1) \ \bigg| \ z'(0) = z'(1) = 0 \right\}$, by $A^\gamma z = (\gamma I - \Delta)z,$
We have that
$\gamma \neq -(n\pi)^2,\quad n = 0, 1, 2, \ldots$. Thus the linear  operator $(A^\gamma)^{-1} : L^2(0, 1) \rightarrow D(A^\gamma)$ is bounded. Then, the solution $(u,v)$ of (\ref{eq:eq1})--(\ref{eq:eq2})  are given by:
\begin{equation}\label{solution v}
    v(\cdot,t)= \beta (\gamma I - \Delta)^{-1}(u(\cdot, t))+ (\gamma I - \Delta)^{-1} f_3(u(\cdot, t)), 
\end{equation}
with  $u$ is given in the invertible transformation (\ref{Ivertible}).\\
\textbf{Step 2: Well-posedness of the parabolic system:}
The parabolic system can be expressed as follows:
	\begin{eqnarray}\label{parabolic}
		\begin{cases}
			\tilde{u_t}(x,t)=\tilde{u}_{xx}(x, t) + \Lambda (\tilde{u}(x, t)),  &   \text{in}\hspace*{0.1cm} (0,1)\times \mathbb{R}^{+},\\
			\tilde{u}_x(0, t) =  \tilde{u}_x(1, t) = 0,  &  \text{in}\hspace*{0.1cm}  \mathbb{R}^{+}, \\
			\tilde{u}(x,0)=\tilde{u}_0, & \text{in}\hspace*{0.1cm} (0,1),
		\end{cases}
	\end{eqnarray} 
 where $\Lambda $ is defined by
\begin{gather*}
\Lambda (\tilde{u}(x, t))= - (c_1 + \rho)\, \tilde{u}(x, t) +\alpha v(x, t) -\\
 \alpha \int_0^x k(x, y)\, v(y, t)\, dy+	\Theta(\tilde{u}(x, t))
 - \int_0^x k(x, y)\Theta(\tilde{u}(y, t))dy,
\end{gather*}
with
\begin{equation}\label{function Theta}
    \Theta(\tilde{u}(x, t))=f_2(v(x, t)) +f_1\left(\tilde{u}(x, t)+\int_0^x  \ell(x, y)\, \tilde{u}(y, t)\, dy\right),
\end{equation}
and 
\[
\begin{aligned}
& v(x,t)= \beta (\gamma I - \Delta)^{-1}\left( \tilde{u}(x, t) + \int_0^x \ell(x, y) \tilde{u}(y, t) \, dy\right)\\
&+ (\gamma I - \Delta)^{-1} f_3\left(\tilde{u}(x, t) + \int_0^x \ell(x, y) \tilde{u}(y, t) \, dy\right),
\end{aligned}
\]
\textbf{i)} \textbf{Show that $\Lambda$ is globally  Lipschitz on $L^2(0,1)$:} Let \( \tilde{u}_1, \tilde{u}_2 \in L^2(0,1) \), and define for each \( x \in (0,1) \),
\[
z_i(x) = \tilde{u}_i(x) + \int_0^x \ell(x, y)\, \tilde{u}_i(y)\, dy, \quad i=1,2.
\]
For all $\tilde{u}_i\in  L^2(0,1)$, we set 
$$
v:=v(\tilde{u}_i)(x) = \beta\, (\gamma I - \Delta)^{-1}(z_i)(x) \\  + (\gamma I - \Delta)^{-1}(f_3(z_i))(x).
$$
We first estimate \( \|z_1 - z_2\| \). Using Lemma \ref{lem:morris2025boundary} and the Cauchy-Schwarz inequality, we get
\begin{gather*}
|z_1(x) - z_2(x)| \leq |\tilde{u}_1(x) - \tilde{u}_2(x)|
+  N_{c_1}\|\tilde{u}_1 - \tilde{u}_2\|.
\end{gather*}
Hence, 
\begin{equation}\label{W1}
    \|z_1 - z_2\| \leq \sqrt{2(1 + N_{c_1}^2)} \, \|\tilde{u}_1 - \tilde{u}_2\|.
\end{equation}
Next, since \( f_3 \) is globally Lipschitz with constant \( M_3 \), and \( (\gamma I - \Delta)^{-1} \) is bounded, we have:
\begin{gather}
\|v(\tilde{u}_1) - v(\tilde{u}_2)\| \leq |\beta|\, \|(\gamma I - \Delta)^{-1}(z_1 - z_2)\| \notag \\
 + \|(\gamma I - \Delta)^{-1}(f_3(z_1) - f_3(z_2))\|\label{W2} \\
\leq \left(|\beta| + M_3\right) \|z_1 - z_2\| \leq L_1 \|\tilde{u}_1 - \tilde{u}_2\|\notag.
\end{gather}
with $L_1=\left(|\beta|  + M_3\right)\sqrt{2(1 + N_{c_1}^2)} $.  \\
Let $\xi(x)=\int_0^x k(x, y)\, (v(\tilde{u}_1(y))\,  -  v(\tilde{u}_2(y))\,) dy$, for all $x\in (0,1)$. Under the Cauchy-Schwarz inequality, we have 
\[
\begin{aligned}
|\xi(x)| &\leq N_{c_1} \|v(\tilde{u}_1) - v(\tilde{u}_2)\|.
\end{aligned}
\]
Using  inequality (\ref{W2}), we get 
$$
\|\xi\| \leq N_{c_1} L_1 \|\tilde{u}_1 - \tilde{u}_2\|.
$$
Furthermore, based on the Assumption \ref{assuption:Lischitz_F}, inequalities (\ref{W1}) and (\ref{W2}), we then conclude that
\[
\begin{aligned}
\|f_1 (z_1) - f_1 (z_2) \|&\leq M_1 \|z_1 - z_2\|\\
&\leq M_1 \sqrt{2(1 + N_{c_1}^2)} \, \|\tilde{u}_1 - \tilde{u}_2\|, 
\end{aligned}
\]
and
\[
\begin{aligned}
\|f_2 (v(\tilde{u}_1)) - f_2 (v(\tilde{u}_2)) \|&\leq M_2 \|v(\tilde{u}_1) - v(\tilde{u}_2)\|\\
&\leq M_2 L_1 \|\tilde{u}_1 - \tilde{u}_2\|.
\end{aligned}
\]
Then 
$$
\|  \Theta(\tilde{u}_1)-  \Theta(\tilde{u}_2)  \|\leq L_2  \|\tilde{u}_1 - \tilde{u}_2\|, 
$$
with $L_2=M_2 L_1 +M_1 \sqrt{2(1 + N_{c_1}^2)}.$
Consequently, we have 
\begin{equation}\label{Lipschitz}
    \|\Lambda(\tilde{u}_1) - \Lambda(\tilde{u}_2)\| \leq L_3\, \|\tilde{u}_1 - \tilde{u}_2\|,
\end{equation}
with  $ L_3 =(1+N_{c_1})(L_1|\alpha|+L_2)+|c_1 + \rho|$. Therefore, the mapping \( \tilde{u} \mapsto \Lambda \) is globally Lipschitz on \( L^2(0,1) \).\\

 \textbf{ii)} \textbf{Existence, uniqueness and regularity of the global mild solution of system (\ref{parabolic}):}\\
Let $A$  be operator defined by $A=\partial_{xx}$ with domaine $D(A) = \left\{ \tilde{u} \in H^2(0, 1) \ \bigg| \ \tilde{u}'(0) = \tilde{u}'(1) = 0 \right\}$. The spectrum of $A$ is given by the	simple eigenvalues $\lambda_{j}= -(j\pi)^{2}$ and eigenfunctions $\varphi_{j}=\sqrt 2 cos(j\pi x)$, $j\geq 0$. Furthermore $A$ generates a contraction semigroup given by $S(t) \tilde{u}=\sum\limits_{j=1}^{\infty} \exp \left(\lambda_{j} t\right)\langle \tilde{u}, \phi_j\rangle \phi_j$, for all $\tilde{u}\in L^{2}(0,1)$ and all $t\geq 0$.\
Therefore, the mapping \( \tilde{u} \mapsto \Lambda \) is globally Lipschitz on \( L^2(0,1) \).
According to  (\cite{pazy}, p. 184), we conclude that for all $\tilde{u}_0\in L^2(0,1)$, the system (\ref{parabolic}) has a unique mild solution
$\tilde{u}\in C([0,+\infty), L^2(0, 1))$, given by 
$$
\tilde{u}(\cdot, t)=S(t)\tilde{u}_0+\int_0^t S(t-\tau)\Lambda (\tilde{u}(\cdot, \tau))d\tau, \quad \forall t \geq 0.
$$
\end{proof}
\subsection{Proof of Theorem~\ref{thm:stability}}
\begin{proof}
We define the Lyapunov function candidate:
		\[
		V(t) = \frac{1}{2} \int_0^1 \tilde{u}^2(x, t)\, dx = \frac{1}{2} \| \tilde{u}(\cdot, t) \|^2.
		\]
        Without loss of generality, we can assume that $V(t)$ is differentiable. Otherwise, we begin with a state in the operator domain for which the convergent solution is differentiable. Finally, we conclude with a density argument.\\
		Taking the time derivative of \( V(t) \), we obtain:
\begin{gather*}
\dot{V}(t) \leq -(c_1 + \rho) \int_0^1 \tilde{u}^2(x, t)\, dx 
		+ \alpha \int_0^1 \tilde{u}(x, t)\, v(x, t)\, dx \\
		- \alpha \int_0^1 \tilde{u}(x, t) \left( \int_0^x k(x, y)\, v(y, t)\, dy \right) dx \\
        + \int_0^1 \Theta(\tilde{u}(x, t)) \tilde{u}(x, t)dx\\
        - \int_0^1 \tilde{u}(x, t) \int_0^x k(x, y) \Theta(\tilde{u}(y, t))dydx,
\end{gather*}
with  $\Theta$ is given in (\ref{function Theta}).      
Using Cauchy–Schwarz inequality, Lemma \ref{lem:morris2025boundary} and Lemma \ref{Lemma5}, we estimate the terms on the right-hand-side of the above inequality  as follows:
$$	\alpha \int_0^1 \tilde{u}(x, t)\, v(x, t)\, dx 
		\leq |\alpha|(|\beta| +M_3 )  (1 +N_{c_1} )\| \tilde{u}(\cdot, t) \|^2,
		$$
	\begin{align*}
		-& \alpha \int_0^1 \tilde{u}(x, t) \left( \int_0^x k(x, y)\, v(y, t)\, dy \right) dx\\
		&\leq |\alpha| \, N_{c_1}  \, \|\tilde{u}(\cdot, t)\| \ \|v(\cdot, t)\| \\
		&\leq |\alpha | N_{c_1} (|\beta| +M_3 )  (1 +N_{c_1} )  \, \|\tilde{u}(\cdot, t)\|^2.
	\end{align*}	
    and
\begin{gather*}
\int_0^1 \Theta(\tilde{u}(x, t)) \tilde{u}(x, t)dx \leq M_2\| v(\cdot, t)) \| \| \tilde{u}(\cdot, t)) \| \\+   M_1 (1+ N_{c_1}) \| \tilde{u}(\cdot, t)  \|^2 
\leq  L_5(1 +N_{c_1} )\| \tilde{u}(\cdot, t)) \|^2,
\end{gather*}
with $L_5=M_1+M_2(|\beta| +M_3 ).  $\\
Hence,  
\begin{align*}
-\int_0^1 \tilde{u}(x, t) &\int_0^x k(x, y) \Theta(\tilde{u}(y, t))dydx \\
& \leq N_{c_1}L_5(1 +N_{c_1} )\|\tilde{u}(\cdot, t)\|^2.
\end{align*}
Then, we deduce that  
		$$
		\dot{V}(t) \leq  -K_1 \|\tilde{u}(\cdot, t)\|^2,
		$$
        with $K_1=c_1 + \rho- ( L_5 +|\alpha | (|\beta| +M_3)) (1 +N_{c_1} )^2 $.
		It follows that
	\begin{equation}\label{stability}
	    \| \tilde{u}(\cdot, t) \| \leq e^{- K_1 t}	\| \tilde{u}_0 \|.
	\end{equation}
    Then 
\[
\begin{aligned}
 \| \tilde{u}(\cdot, t) \| + \| v(\cdot, t) \| 
&\leq \sqrt{2} \left\|
\begin{bmatrix}
\tilde{w}(\cdot, t) \\
v(\cdot, t)
\end{bmatrix}
\right\|_{L^2 \times L^2}\\
&\leq \sqrt{2} \big(1+ (N_{c_1}+1)(|\beta| +M_3 )\big) \| \tilde{u}(\cdot, t) \|.
\end{aligned}
\]
By inequality (\ref{stability}) and Lemma~\ref{Lemma5}, we find that
	 \[
	 \| \tilde{u}(\cdot, t) \| + \| v(\cdot, t) \| \leq  c \left( \| \tilde{u}_0 \| + \| v_0 \| \right) e^{-K_1 t} ,
	 \]
     with $c=\sqrt{2} (1+ (N_{c_1}+1)(|\beta| +M_3 )).$
        \end{proof}
\subsection{ Proof of Theorem \ref{Theo of sys1}}
\begin{proof}
    According to  Theorem \ref{thm:stability} that the target system (\ref{eq:eq755}) is exponentially stable.  Furthermore, by Theorem \ref{thm:target_system} there is an invertible state transformation between the original system (\ref{eq:eq1})--(\ref{eq:eq4}) with $w(t)$ given by (\ref{control}) and the exponentially stable target system (\ref{eq:eq755}). \\
Note that the kernel functions \( k(x, y) \) and \( \ell(x, y) \) are bounded. Then, a straightforward generalization of (Theorem 4, \cite{SmyshlyaevKrstic2004})~ yields that there is an equivalence of norms between \( u(., t) \) and \( \tilde{u}(., t) \) in \( L^2(0,1) \). Then, we can derive similar convergence properties for the original system (\ref{eq:eq1})--(\ref{eq:eq4}) as the ones for the target system since the backstepping transformation is invertible.   This completes the proof.
\end{proof}

\subsection{Proof of Proposition \ref{obsever_pro}}
\begin{proof}
We first take the spatial derivatives of  (\ref{tran11}):
\begin{gather}
\er^{u}_{xx}(x, t) = \er^{\tilde{u}}_{xx}(x, t) + \int_0^x k_{xx}(x, y)\, \er^{\tilde{w}}(y, t)\, dy 
 \\ + k^a_x(x, x)\, \er^{\tilde{u}}(x, t) + \frac{d}{dx}k(x, x)\, \er^{\tilde{u}}(x, t) 
 + k(x, x)\, \er^{\tilde{u}}_x(x, t).\notag
\label{dérivé_xx}
\end{gather}
Taking the time derivative of (\ref{tran11}) and integrating by parts, we get:
\begin{align}
 &\er^{u}_t(x, t)= \er^{\tilde{u}}_t(x, t) + \int_0^x k(x, y)\, \er^{\tilde{u}}_t(y, t)\, dy \notag \\
&= \er^{\tilde{u}}_t(x, t) - \rho \int_0^x k(x, y)\, \er^{\tilde{u}}(y, t)\, dy \notag  \\
&\quad  + \alpha \int_0^x k(x, y)\, \er^{v}(y, t)\, dy + k(x, x)\, \er^{\tilde{u}}_x(x, t) \notag  \\
&\quad - k(x, 0)\, \er^{\tilde{u}}_x(0, t)- k_y(x, x)\, \er^{\tilde{u}}(x, t) + k_y(x, 0)\, \er^{\tilde{u}}(0, t) \notag \\
&\quad + \int_0^x k_{yy}(x, y)\, \er^{\tilde{u}}(y, t)\, dy - \er^{\tilde{u}}(1, t)\int_0^x k(x, y)\, \sigma_1(y)\, dy\notag \\
&\quad  +  \int_0^x k(x, y)  \left( f_1(u(y, t))  -f_1(u(y, t)-\er^{u}(y, t) ) \right) dy \notag \\
&\quad   +\int_0^x k(x, y)  \left( f_2(v(y, t))  -f_2(v(y, t)-\er^{v}(y, t) ) \right)dy\notag \\
&\quad  + f_1(u(x, t))  -f_1(u(x, t)-\er^{u}(x, t) )  \notag \\
&\quad  +  f_2(v(x, t))  -f_2(v(x, t)-\er^{v}(x, t)). 
\label{dérivé_t}
\end{align}
Substituting (\ref{dérivé_xx}) and (\ref{dérivé_t}) in the parabolic equation (12), and using $\er^{\tilde{u}}_x(0, t)=0$  and $ k_y(x, 0)=0$, we obtain that 
{\small
\begin{gather}
\er_t^{\tilde{w}}(x, t) = \er_{xx}^{\tilde{w}}(x, t) + \alpha\, \er^v(x, t) + \left(k_y^a(x, x) + k_x^a(x, x)\right) \er^u(x, t) \notag \\
 + k(x, x)\, \er_x^u(x, t) - k(x, x)\, \er_x^u(x, t) - \alpha \int_0^x k(x, y)\, \er^v(y, t)\, dy \notag \\
 - \sigma_1(x)\, \er^u(1, t) + \int_0^x \left[k_{xx}(x, y) - k_{yy}(x, y) + \rho\, k(x, y)\right] \er^u(y, t)\, dy \notag \\
+ \er^u(1, t) \int_0^x k(x, y)\, \sigma_1(y)\, dy + \left(-\rho + \frac{d}{dx}k(x, x)\right) e(x, t) \notag \\
 - \int_0^x k(x, y) \left( f_1(u(y, t)) - f_1(u(y, t) - \er^u(y, t)) \right)\, dy \notag \\
 - \int_0^x k(x, y) \left( f_2(v(y, t)) - f_2(v(y, t) - \er^v(y, t)) \right)\, dy \notag \\
- f_1(u(x, t)) + f_1(u(x, t) - \er^u(x, t)) \notag \\
 - f_2(v(x, t)) + f_2(v(x, t) - \er^v(x, t)).
\label{eq:etildewn}
\end{gather}}
We add and subtract the term $c_1 \er^{\tilde{w}}$ to the right-hand-side of equation (\ref{eq:etildewn}). Moreover, we have that  $k(x, y)$ is given by  Lemma \ref{lem:SmyshlyaevKrstic2004}  and  $\sigma_1(x)=0$, then  the previous equation yields
\begin{gather}
\er_t^{\tilde{u}}(x, t) = \er_{xx}^{\tilde{u}}(x, t) -(c_1+\rho) \er^{\tilde{u}}(x, t) + \alpha\, \er^v(x, t)  - \notag \\
 \alpha \int_0^x k(x, y)\, \er^v(y, t)\, dy -f_2(v(x, t)) + f_2(v(x, t) - \er^v(x, t))\notag  \\
 - \int_0^x k(x, y) \left( f_1(u(y, t)) - f_1(u(y, t) - \er^u(y, t)) \right)\, dy \notag \\
 - \int_0^x k(x, y) \left( f_2(v(y, t)) - f_2(v(y, t) - \er^v(y, t)) \right)\, dy \notag \\
 - f_1(u(x, t)) + f_1(u(x, t) - \er^u(x, t)).
\label{eq:etildew}
\end{gather}
 Furthermore, from Lemma \ref{lem:SmyshlyaevKrstic2004}, we get 
$$
\er_x^{\tilde{u}}(0, t)=\er_x^u(0, t)-\int_0^0 k_x(1, y) u(y,t) d y-k(0,0) u(1, t)=0,
$$
\begin{gather}
 \er_x^{\tilde{u}}(1, t) =\er_x^u(1, t)-\int_0^1 k_x(1, y) \er^u(y, t) d y-k(1,1) \er^u(1, t)\notag \\
 =-\left(\sigma_2+k(1,1)\right) \er^w(1, t)-\int_0^1 k_x(1, y) \er^u(y, t) d y \notag\\
 =-\int_0^1 k_x(1, y) \er^u(y, t) d y 
 =-\int_0^1 k_x(1, y) \er^{\tilde{u}}(y, t) d y\notag\\-\int_0^1 k_x(1, y) \int_0^y \ell(y, z) \er^{\tilde{u}}(z, t) d z d y
\end{gather}
\end{proof}
\vspace{-0.5cm}
\subsection{Proof of Theorem \ref{stability_Observer}}
With a parallel line of reasoning as the one used to prove Theorem~\ref{thm:stability}, it follows that the stability of the target system (\ref{eq:o11})-(\ref{eq:o3}) follows from the exponential decay of the state $\er^{\tilde{u}}(x, t)$. Define the Lyapunov function candidate
\[
		V(t) = \int_0^1 (\er^{\tilde{u}}(x, t))^2 dx = \frac{1}{2} \| \er^{\tilde{u}}(\cdot, t) \|^2.
		\]
		Taking the time derivative of \( V(t) \), we obtain:
\begin{gather*}
\dot{V}(t)  =\int_0^1 \er^{\tilde{u}}(x, t) \er_{x x}^{\tilde{u}}(x, t) d x-\left(c_1+\rho\right) \int_0^1\left(\er^{\tilde{u}}(x, t)\right)^2 d x \\
 +\alpha \int_0^1 \er^{\tilde{u}}(x, t) \er^v(x, t) d x-\alpha \int_0^1 \er^{\tilde{u}}(x, t) \int_0^x k(x, y)   \er^v(y, t)\\ d y d x
  -\int_0^1 \er^{\tilde{u}}(x, t)\int_0^x k(x, y) G_1(\er^{u}(y,t))dydx\\
 -\int_0^1 \er^{\tilde{u}}(x, t)\left( G_1(\er^{u}(x,t))+G_2(\er^{v}(x,t))\right)dx\\
-\int_0^1 \er^{\tilde{u}}(x, t)\int_0^x k(x, y) G_2(\er^{v}(y,t))dydx,
\end{gather*}
where$G_1$ and $G_2$ are defined as follows:
$$
G_1(\er^{u}(x, t) )= f_1(u(x, t))  -f_1(u(x, t)-\er^{u}(x, t) ),
$$
$$
G_2(\er^{v}(x, t))= f_2(v(x, t))  -f_2(v(x, t)-\er^{v}(x, t) ).
$$
 We have that $\gamma \neq -(n\pi)^2 $ \text{with}  $n\in \mathbb{N}$, then 
\[
\begin{aligned}
 \er^{v}(\cdot,t)=&  (\gamma I - \Delta)^{-1} \bigg(f_3(u(\cdot, t))-f_3(u(\cdot, t)-\er^{u}(\cdot, t)) \bigg)\\
 &+\beta (\gamma I - \Delta)^{-1} (\er^{u}(\cdot, t)).    
 \end{aligned}\]
The norms of the functions $v$, $G_1$ and $G_2$ in $L^2$ are given as follows:

\begin{equation}\label{46}
\begin{aligned}
 \| \er^{v}(\cdot,t) \| &\leq | \beta |\|  \er^{u}(\cdot, t) \|+\|  f_3(u(\cdot, t))-f_3(u(\cdot, t)-\er^{u}(\cdot, t))\|\\
&\leq (| \beta |+M_3)\|  \er^{u}(\cdot, t) \|\\
 &\leq (| \beta |+M_3)(1+N_{c_1})\| \er^{\tilde{u}}(\cdot, t) \|
\end{aligned}
\end{equation}
 and
 \begin{equation}\label{47}
\begin{aligned}
 \| G_1(\er^{u}(\cdot, t) ) \| &\leq M_1\| \er^{u}(\cdot, t) \|\\
 &\leq M_1(1+N_{c_1})\| \er^{\tilde{u}}(\cdot, t) \|.
 \end{aligned}
 \end{equation}
 
From (\ref{46}), we deduce that 
\begin{equation}\label{48}
\begin{aligned}
 \| G_2(\er^{v}(\cdot, t) ) \| &\leq M_2\| \er^{v}(\cdot, t) \|\\
 &\leq  M_2(| \beta |+M_3)(1+N_{c_1})\| \er^{\tilde{u}}(\cdot, t) \|.
 \end{aligned}
\end{equation}
Using the  Cauchy–Schwarz inequality and (\ref{46}), we have
$$	\alpha \int_0^1 \er^{\tilde{u}}(x, t)\, \er^v(x, t)\, dx 
		\leq |\alpha|(|\beta| +M_3 )  (1 +N_{c_1} )\| \er^{\tilde{u}}(\cdot, t) \|^2,
		$$
	 \begin{align*} \text{and} \ 
		&- \alpha \int_0^1 \er^{\tilde{u}}(x, t) \left( \int_0^x k(x, y)\, \er^{v}(y, t)\, dy \right) dx\\
		&\leq |\alpha| \, N_{c_1}  \, \|\er^{\tilde{u}}(\cdot, t)\| \, \|\er^{\tilde{v}}(\cdot, t)\| \\
		&\leq |\alpha | N_{c_1} (|\beta| +M_3 )  (1 +N_{c_1} )  \, \|\er^{\tilde{u}}(\cdot, t)\|^2.
	\end{align*}	
From (\cite{morris2025boundary}, p. 3), we have that 
\[
\begin{aligned}
 \int_0^1 \er^{\tilde{u}}(x, t) \er_{x x}^{\tilde{u}}(x, t) d x  \leq \frac{(1+N_{c_1})^2\eta^2+1}{2}\left\|\er^{\tilde{u}}(\cdot, t)\right\|^2,
\end{aligned}\]
with $\eta= \frac{c_1}{2}\left(1+\frac{c_1}{2}\right) e^{\frac{c_1}{4}}\left(\sqrt{\frac{\pi}{2 c_1}} \operatorname{erf}\left(\sqrt{\frac{c_1}{2}}\right)\right)^{\frac{1}{2}}$.\\
Using  (\ref{47})  and (\ref{48}), we get that, via  Cauchy–Schwarz inequality 
\[\begin{aligned}
 &-\int_0^1 \er^{\tilde{u}}(x, t)\int_0^x k(x, y) G_1(\er^{u}(y,t))dydx\\ &\leq N_{c_1}\| \er^{\tilde{u}} (\cdot, t) \| \| G_1(\er^{u} (\cdot, t)) \|\\
 &\leq N_{c_1} M_1(1+N_{c_1})\| \er^{\tilde{u}}(\cdot, t) \|^2,
 \end{aligned} \]
 and 
\begin{equation*}
\begin{aligned}
 &-\int_0^1 \er^{\tilde{u}}(x, t)\int_0^x k(x, y) G_2(\er^{v}(y,t))dydx\\ &\leq N_{c_1}\| \er^{\tilde{u}} (\cdot, t) \| \| G_2(\er^{v} (\cdot, t)) \|\\
 &\leq N_{c_1} M_2(| \beta |+M_3)(1+N_{c_1})\| \er^{\tilde{u}}(\cdot, t) \|^2,
 \end{aligned}
 \end{equation*}
 and that 
 \begin{equation*}\label{50}
\begin{aligned}
 &-\int_0^1 \er^{\tilde{u}}(x, t)\left( G_1(\er^{u}(x,t))+G_2(\er^{v}(x,t))\right)dx\\ &\leq \| \er^{\tilde{u}} (\cdot, t) \| \| G_2(\er^{v} (\cdot, t))+G_1(\er^{u} (\cdot, t)) \|\\
 &\leq \big(M_2(| \beta |+M_3)+M_1\big) (1+N_{c_1})\| \er^{\tilde{u}}(\cdot, t) \|^2,
 \end{aligned}
\end{equation*}
Using the inequalities previously mentioned leads to the following result:
$$
\dot{V}(t)\leq -2 K_3V(t), \ \text{with} \ 
k_3=c_1+\rho -K_4,
$$

and 
$$
 K_4= \left(( M_2 +|\alpha |  ) (|\beta| +M_3 ) +\frac{\eta^2+1}{2}+ M_1\right) (1+N_{c_1})^2.
$$
Similarly to the   Theorem \ref{thm:stability},   if the parameter $c_1$ is chosen such that (\ref{condition}) is satisfied, then the observation error system (\ref{eq:obs_err1})-(\ref{eq:obs_errn})  is exponentially stable.

\section{Conclusions}
In this paper, we have developed a boundary control and state observer design for a class of nonlinear coupled parabolic-elliptic systems. The proposed approach is based on the backstepping technique. Several conditions are provided to ensure well-posedness and exponential stability in the presence of nonlinearities. In future work, we aim to relax the assumptions imposed on the nonlinear functions.
\bibliographystyle{IEEEtran}
\balance
\bibliography{biblios}

\begin{thebibliography}{10}
\providecommand{\url}[1]{#1}
\csname url@samestyle\endcsname
\providecommand{\newblock}{\relax}
\providecommand{\bibinfo}[2]{#2}
\providecommand{\BIBentrySTDinterwordspacing}{\spaceskip=0pt\relax}
\providecommand{\BIBentryALTinterwordstretchfactor}{4}
\providecommand{\BIBentryALTinterwordspacing}{\spaceskip=\fontdimen2\font plus
\BIBentryALTinterwordstretchfactor\fontdimen3\font minus
  \fontdimen4\font\relax}
\providecommand{\BIBforeignlanguage}[2]{{%
\expandafter\ifx\csname l@#1\endcsname\relax
\typeout{** WARNING: IEEEtran.bst: No hyphenation pattern has been}%
\typeout{** loaded for the language `#1'. Using the pattern for}%
\typeout{** the default language instead.}%
\else
\language=\csname l@#1\endcsname
\fi
#2}}
\providecommand{\BIBdecl}{\relax}
\BIBdecl

\bibitem{Ahn2019}
J.~Ahn and C.~Yoon, ``Global well-posedness and stability of constant
  equilibria in parabolic–elliptic chemotaxis systems without gradient
  sensing,'' \emph{Nonlinearity}, vol.~32, no.~4, p. 1327, 2019.

\bibitem{Wu2006}
J.~Wu, J.~Xu, and H.~Zou, ``On the well-posedness of a mathematical model for
  lithium-ion battery systems,'' \emph{Methods and Applications of Analysis},
  vol.~13, no.~3, pp. 275--298, 2006.

\bibitem{TaoWinkler2015}
Y.~Tao and M.~Winkler, ``Boundedness vs. blow-up in a two-species chemotaxis
  system with two chemicals,'' \emph{Discrete and Continuous Dynamical Systems
  - B}, vol.~20, no.~9, pp. 3165--3183, 2015.

\bibitem{Parada2020}
H.~Parada, E.~Cerpa, and K.~Morris, ``Feedback control of an unstable
  parabolic-elliptic system with input delay,'' 2020, preprint.

\bibitem{morris2025boundary}
A.~Alalabi and K.~Morris, ``Boundary control and observer design via
  backstepping for a coupled parabolic--elliptic system,'' \emph{Automatica},
  vol. 174, p. 112154, 2025.

\bibitem{KrsticSmyshlyaev2008}
M.~Krstic and A.~Smyshlyaev, \emph{{Boundary Control of {PDEs}: A Course on
  Backstepping Designs}}.\hskip 1em plus 0.5em minus 0.4em\relax SIAM, 2008.

\bibitem{VazquezKrstic2009}
R.~Vazquez and M.~Krstic, ``Boundary observer for output-feedback stabilization
  of thermal-fluid convection loop,'' \emph{IEEE Transactions on Control
  Systems Technology}, vol.~18, no.~4, pp. 789--797, 2009.

\bibitem{VazquezTrelatCoron2008}
R.~V. Valenzuela, E.~Trélat, and J.~M. Coron, ``Control for fast and stable
  laminar-to-high-reynolds-numbers transfer in a {2D Navier–Stokes} channel
  flow,'' \emph{Discrete and Continuous Dynamical Systems - B}, vol.~10, no.~4,
  pp. 925--956, 2008.

\bibitem{VazquezKrstic2006}
R.~Vazquez and M.~Krstic, ``Explicit integral operator feedback for local
  stabilization of nonlinear thermal convection loop {PDEs},'' \emph{Systems \&
  Control Letters}, vol.~55, no.~8, pp. 624--632, 2006.

\bibitem{pazy}
A.~Pazy, \emph{Semigroups of linear operators and applications to partial
  differential equations}.\hskip 1em plus 0.5em minus 0.4em\relax Springer
  Science \& Business Media, 2012, vol.~44.

\bibitem{SmyshlyaevKrstic2004}
A.~Smyshlyaev and M.~Krstic, ``Closed-form boundary state feedbacks for a class
  of {1-D} partial integro-differential equations,'' \emph{IEEE Transactions on
  Automatic Control}, vol.~49, no.~12, pp. 2185--2202, 2004.

\end{thebibliography}
\appendix 
\subsection{Useful results} \label{appendixA}
 In this section, we recall and prove some results that are used in the paper, which provide useful estimates for the stabilization problem.
\begin{lemma}(\cite{SmyshlyaevKrstic2004}, chap. 4)\label{lem:SmyshlyaevKrstic2004}
		For any $c_1 > 0$, the following boundary value problem:
		\begin{subequations}
			\[
			\begin{aligned}
				&k_{yy}(x, y) - k_{xx}(x, y) + c_1 k(x, y) = 0, & 0 < y < x < 1,  \\
				&k_y(x, 0) = 0, \quad k(x, x) = -\frac{1}{2}c_1 x, 
			\end{aligned}
            \]
		\end{subequations}
		has a continuous unique solution.
        \end{lemma}
\begin{lemma}(\cite{SmyshlyaevKrstic2004}, chap. 4)\label{Lemma3}  Assuming that $c_1> 0$, the inverse transformation of (\ref{KAMAL}) is
\begin{equation}\label{Ivertible}
u(x, t) = \tilde{u}(x, t) + \int_0^x \ell(x, y) \tilde{u}(y, t) \, dy,
\end{equation}
where$\ell(x, y)$ is the solution of the system
\begin{align}
\ell_{xx}(x, y) - \ell_{yy}(x, y) + c_1^2 \ell(x, y) &= 0, \\
\ell_{y}(x, 0) = 0, \quad \ell(x, x) &= -\frac{1}{2}c_1^2 x. 
\end{align}
\end{lemma}
        
\begin{lemma}(\cite{morris2025boundary}, p. 3)\label{lem:morris2025boundary}
  For  $c_1> 0$, letting $$N_{c_1} = \sqrt{\frac{c_1 \pi}{8}} \left( \operatorname{erfi}(\sqrt{\frac{2}{c_1}}) \cdot \operatorname{erf}(\sqrt{\frac{2}{c_1}}) \right)^{\frac{1}{2}},$$
				wherethe functions $\operatorname{erfi}(x)$ and $\operatorname{erf}(x)$ are defined as:
				\[
				\operatorname{erfi}(x) = \frac{2}{\sqrt{\pi}} \int_0^x e^{s^2}\, ds, \qquad 
				\operatorname{erf}(x) = \frac{2}{\sqrt{\pi}} \int_0^x e^{-s^2}\, ds.
				\]
				
				The $L^2$-norms of the kernels $k$ and $\ell$ are bounded by:
				\[
				\|k\|\leq N_{c_1}, \qquad \|\ell\|\leq N_{c_1},
				\]
                \end{lemma}
The following lemma, which establishes a relation of the $L^2$-norm of $v(\cdot, t)$ to that of $\tilde{u}(\cdot, t)$, will be needed to show the exponential  stability of the target system.
\begin{lemma}\label{Lemma5}
    Let $\gamma \neq -(n\pi)^2 $ \text{with}  $n\in \mathbb{N}$. The solution of (\ref{eq:eq755}) satisfies
\begin{equation}\label{norme v}
    \|v(\cdot, t)\|\leq (|\beta| +M_3 )  (1 +N_{c_1} ) \|\tilde{u}(\cdot, t)\|. 
\end{equation}
\end{lemma}
\begin{remark}
    In the above lemma, we established a relation between the norms of the solutions in the general case where$\gamma \neq -(n\pi)^2 $, without requiring  $\gamma$ to be strictly positive. This extends the result obtained by \cite{morris2025boundary}, who proved a similar norm relation under the more restrictive assumption  $\gamma >0$.
\end{remark}
\begin{lemma}\label{error_sys_sol}
Let $\gamma \neq -(n\pi)^2 $ and assume Assumption \ref{assuption:Lischitz_F} holds. Then, the system error observer (12)-(15) admits  unique mild solution $(\er^{u}, \er^{v})\in C([0,+\infty), L^2(0, 1)) \times C([0,+\infty), L^2(0, 1))$.
\end{lemma}
\vspace{-0.5cm}
\subsection{Proof of lemmas}\label{appendixB}
\subsubsection{Proof of Lemma~\ref{Lemma5}}
\begin{proof}
According to Theorem~\ref{thm:well-posdnes}, the initial system \eqref{eq:eq1}--\eqref{eq:eq4} possesses a unique solution such that:
$$
v(\cdot,t)= \beta (\gamma I - \Delta)^{-1}(u(\cdot, t))+ (\gamma I - \Delta)^{-1} f_3(u(\cdot, t)), 
$$
with  $u$ is given in (\ref{Ivertible}).
By Cauchy–Schwarz inequality  and  Lemma \ref{lem:morris2025boundary}, we deduce that  
\begin{equation}\label{U1}
\|u(\cdot, t)\|\leq (1+N_{c_1}) \|\tilde{u}(\cdot, t)\|.
\end{equation}
Using the fact that $(\gamma I - \Delta)^{-1}$ is bounded, under Assumption (\ref{assuption:Lischitz_F})  and (\ref{U1}), we can get~\eqref{norme v}.
 \end{proof}
 \subsubsection{Proof of Lemma~\ref{error_sys_sol}}
\begin{proof}
Let $n\in \mathbb{N}$,  for $\gamma \neq -(n\pi)^2 $, we have that 
$$
\er^{v}(\cdot,t)= \beta (\gamma I - \Delta)^{-1}(\er^{u}(\cdot, t))+ (\gamma I - \Delta)^{-1} f_3(\er^{u}(\cdot, t)). 
$$
    We define  the operator $\mathcal{A}$ by 
 \begin{equation*}
     \mathcal{A} = \Delta - \rho I,
 \end{equation*}
	with domain
	$$
	D(\mathcal{A}) = \left\{ w \in H^2(0, 1) \ \big| \ w'(0)  = 0 , \quad w'(1)=-\eta w(1)\right\},
	$$      
where$\eta\geq 0$.  The operator $\mathcal{A}$  generates the $C_0$-semigroup $S(t)$ on $L^2(0,1)$.  Furthermore,  we have that $f_1$, $f_2$ and  $f_3$ 
 are globally Lipschitz in  $L^2(0,1)$, then the maps defined by 
 $ F(z)= f_1(u)  -f_1(u-z)+f_2(v)-f_2(v-z)$ is globally Lipschitz in  $L^2(0,1)$, which leads to the conclusion that system (12)-(15) is well posed  ( see \cite{pazy}, p. 184).\\
\end{proof}
 
\end{document}